\documentclass[12pt]{amsart}
\usepackage{amsmath,amssymb,amsfonts}
\begin{document}
\theoremstyle{plain}
\newtheorem{thm}{Theorem}[subsection]
\newtheorem{lem}[thm]{Lemma}
\newtheorem{cor}[thm]{Corollary}
\newtheorem{prop}[thm]{Proposition}

\theoremstyle{definition}
\newtheorem{rem}[thm]{Remark}
\newtheorem{defn}[thm]{Definition}
\newtheorem{ex}[thm]{Example}

\numberwithin{equation}{subsection}
\newcommand{\mc}{\mathcal}
\newcommand{\mb}{\mathbb}
\newcommand{\surj}{\twoheadrightarrow}
\newcommand{\inj}{\hookrightarrow}
\newcommand{\zar}{{\rm zar}}
\newcommand{\an}{{\rm an}} 
\newcommand{\red}{{\rm red}}
\newcommand{\codim}{{\rm codim}}
\newcommand{\rank}{{\rm rank}}
\newcommand{\Ker}{{\rm Ker \ }}
\newcommand{\Pic}{{\rm Pic}}
\newcommand{\Div}{{\rm Div}}
\newcommand{\Hom}{{\rm Hom}}
\newcommand{\im}{{\rm im}}
\newcommand{\Spec}{{\rm Spec \,}}
\newcommand{\Sing}{{\rm Sing}}
\newcommand{\Char}{{\rm char}}
\newcommand{\Tr}{{\rm Tr}}
\newcommand{\Gal}{{\rm Gal}}
\newcommand{\Min}{{\rm Min \ }}
\newcommand{\Max}{{\rm Max \ }}
\newcommand{\sA}{{\mathcal A}}
\newcommand{\sB}{{\mathcal B}}
\newcommand{\sC}{{\mathcal C}}
\newcommand{\sD}{{\mathcal D}}
\newcommand{\sE}{{\mathcal E}}
\newcommand{\sF}{{\mathcal F}}
\newcommand{\sG}{{\mathcal G}}
\newcommand{\sH}{{\mathcal H}}
\newcommand{\sI}{{\mathcal I}}
\newcommand{\sJ}{{\mathcal J}}
\newcommand{\sK}{{\mathcal K}}
\newcommand{\sL}{{\mathcal L}}
\newcommand{\sM}{{\mathcal M}}
\newcommand{\sN}{{\mathcal N}}
\newcommand{\sO}{{\mathcal O}}
\newcommand{\sP}{{\mathcal P}}
\newcommand{\sQ}{{\mathcal Q}}
\newcommand{\sR}{{\mathcal R}}
\newcommand{\sS}{{\mathcal S}}
\newcommand{\sT}{{\mathcal T}}
\newcommand{\sU}{{\mathcal U}}
\newcommand{\sV}{{\mathcal V}}
\newcommand{\sW}{{\mathcal W}}
\newcommand{\sX}{{\mathcal X}}
\newcommand{\sY}{{\mathcal Y}}
\newcommand{\sZ}{{\mathcal Z}}
\newcommand{\A}{{\mathbb A}}
\newcommand{\B}{{\mathbb B}}
\newcommand{\C}{{\mathbb C}}
\newcommand{\D}{{\mathbb D}}
\newcommand{\E}{{\mathbb E}}
\newcommand{\F}{{\mathbb F}}
\newcommand{\G}{{\mathbb G}}
\renewcommand{\H}{{\mathbb H}}
\newcommand{\I}{{\mathbb I}}
\newcommand{\J}{{\mathbb J}}
\newcommand{\M}{{\mathbb M}}
\newcommand{\N}{{\mathbb N}}
\renewcommand{\P}{{\mathbb P}}
\newcommand{\Q}{{\mathbb Q}}
\newcommand{\R}{{\mathbb R}}
\newcommand{\T}{{\mathbb T}}
\newcommand{\U}{{\mathbb U}}
\newcommand{\V}{{\mathbb V}}
\newcommand{\W}{{\mathbb W}}
\newcommand{\X}{{\mathbb X}}
\newcommand{\Y}{{\mathbb Y}}
\newcommand{\Z}{{\mathbb Z}}

\newcommand{\ffa}{{\mathfrak a}}
\newcommand{\ffb}{{\mathfrak b}}
\newcommand{\ffc}{{\mathfrak c}}
\newcommand{\ffd}{{\mathfrak d}}
\newcommand{\ffe}{{\mathfrak e}}
\newcommand{\fff}{{\mathfrak f}}
\newcommand{\ffg}{{\mathfrak g}}
\newcommand{\ffh}{{\mathfrak h}}
\newcommand{\ffi}{{\mathfrak i}}
\newcommand{\ffj}{{\mathfrak j}}
\newcommand{\ffk}{{\mathfrak k}}
\newcommand{\ffl}{{\mathfrak l}}
\newcommand{\ffm}{{\mathfrak m}}
\newcommand{\ffn}{{\mathfrak n}}
\newcommand{\ffo}{{\mathfrak o}}
\newcommand{\ffp}{{\mathfrak p}}
\newcommand{\ffq}{{\mathfrak q}}
\newcommand{\ffr}{{\mathfrak r}}
\newcommand{\ffs}{{\mathfrak s}}
\newcommand{\fft}{{\mathfrak t}}
\newcommand{\ffu}{{\mathfrak u}}
\newcommand{\ffv}{{\mathfrak v}}
\newcommand{\ffw}{{\mathfrak w}}
\newcommand{\ffx}{{\mathfrak x}}
\newcommand{\ffy}{{\mathfrak y}}
\newcommand{\ffz}{{\mathfrak z}}
\newcommand{\ffA}{{\mathfrak A}}
\newcommand{\ffB}{{\mathfrak B}}
\newcommand{\ffC}{{\mathfrak C}}
\newcommand{\ffD}{{\mathfrak D}}
\newcommand{\ffE}{{\mathfrak E}}
\newcommand{\ffF}{{\mathfrak F}}
\newcommand{\ffG}{{\mathfrak G}}
\newcommand{\ffH}{{\mathfrak H}}
\newcommand{\ffI}{{\mathfrak I}}
\newcommand{\ffJ}{{\mathfrak J}}
\newcommand{\ffK}{{\mathfrak K}}
\newcommand{\ffL}{{\mathfrak L}}
\newcommand{\ffM}{{\mathfrak M}}
\newcommand{\ffN}{{\mathfrak N}}
\newcommand{\ffO}{{\mathfrak O}}
\newcommand{\ffP}{{\mathfrak P}}
\newcommand{\ffQ}{{\mathfrak Q}}
\newcommand{\ffR}{{\mathfrak R}}
\newcommand{\ffS}{{\mathfrak S}}
\newcommand{\ffT}{{\mathfrak T}}
\newcommand{\ffU}{{\mathfrak U}}
\newcommand{\ffV}{{\mathfrak V}}
\newcommand{\ffW}{{\mathfrak W}}
\newcommand{\ffX}{{\mathfrak X}}
\newcommand{\ffY}{{\mathfrak Y}}
\newcommand{\ffZ}{{\mathfrak Z}}


\catcode`\@=11
%
%
\def\opn#1#2{\def#1{\mathop{\kern0pt\fam0#2}\nolimits}} 
\def\bold#1{{\bf #1}}%
\def\underrightarrow{\mathpalette\underrightarrow@}
\def\underrightarrow@#1#2{\vtop{\ialign{$##$\cr
 \hfil#1#2\hfil\cr\noalign{\nointerlineskip}%
 #1{-}\mkern-6mu\cleaders\hbox{$#1\mkern-2mu{-}\mkern-2mu$}\hfill
 \mkern-6mu{\to}\cr}}}
\let\underarrow\underrightarrow
\def\underleftarrow{\mathpalette\underleftarrow@}
\def\underleftarrow@#1#2{\vtop{\ialign{$##$\cr
 \hfil#1#2\hfil\cr\noalign{\nointerlineskip}#1{\leftarrow}\mkern-6mu
 \cleaders\hbox{$#1\mkern-2mu{-}\mkern-2mu$}\hfill
 \mkern-6mu{-}\cr}}}
%
%

%
\def\:{\colon}
\let\oldtilde=\tilde
\def\tilde#1{\mathchoice{\widetilde{#1}}{\widetilde{#1}}%
{\indextil{#1}}{\oldtilde{#1}}}
\def\indextil#1{\lower2pt\hbox{$\textstyle{\oldtilde{\raise2pt%
\hbox{$\scriptstyle{#1}$}}}$}}
\def\pnt{{\raise1.1pt\hbox{$\textstyle.$}}}
%

%
\let\amp@rs@nd@\relax
\newdimen\ex@
\ex@.2326ex
\newdimen\bigaw@
\newdimen\minaw@
\minaw@16.08739\ex@
\newdimen\minCDaw@
\minCDaw@2.5pc
\newif\ifCD@
\def\minCDarrowwidth#1{\minCDaw@#1}
\newenvironment{CD}{\@CD}{\@endCD}
\def\@CD{\def\A##1A##2A{\llap{$\vcenter{\hbox
 {$\scriptstyle##1$}}$}\Big\uparrow\rlap{$\vcenter{\hbox{%
$\scriptstyle##2$}}$}&&}%
\def\V##1V##2V{\llap{$\vcenter{\hbox
 {$\scriptstyle##1$}}$}\Big\downarrow\rlap{$\vcenter{\hbox{%
$\scriptstyle##2$}}$}&&}%
\def\={&\hskip.5em\mathrel
 {\vbox{\hrule width\minCDaw@\vskip3\ex@\hrule width
 \minCDaw@}}\hskip.5em&}%
\def\verteq{\Big\Vert&&}%
\def\noarr{&&}%
\def\vspace##1{\noalign{\vskip##1\relax}}\relax\let\amp@rs@nd@&\iffalse}\fi
 \CD@true\vcenter\bgroup\relax\let\\=\cr\iffalse}\fi\tabskip\z@skip\baselineskip20\ex@
 \lineskip3\ex@\lineskiplimit3\ex@\halign\bgroup
 &\hfill$\m@th##$\hfill\cr}
\def\@endCD{\cr\egroup\egroup}
%
\def\>#1>#2>{\amp@rs@nd@\setbox\z@\hbox{$\scriptstyle
 \;{#1}\;\;$}\setbox\@ne\hbox{$\scriptstyle\;{#2}\;\;$}\setbox\tw@
 \hbox{$#2$}\ifCD@
 \global\bigaw@\minCDaw@\else\global\bigaw@\minaw@\fi
 \ifdim\wd\z@>\bigaw@\global\bigaw@\wd\z@\fi
 \ifdim\wd\@ne>\bigaw@\global\bigaw@\wd\@ne\fi
 \ifCD@\hskip.5em\fi
 \ifdim\wd\tw@>\z@
 \mathrel{\mathop{\hbox to\bigaw@{\rightarrowfill}}\limits^{#1}_{#2}}\else
 \mathrel{\mathop{\hbox to\bigaw@{\rightarrowfill}}\limits^{#1}}\fi
 \ifCD@\hskip.5em\fi\amp@rs@nd@}
\def\<#1<#2<{\amp@rs@nd@\setbox\z@\hbox{$\scriptstyle
 \;\;{#1}\;$}\setbox\@ne\hbox{$\scriptstyle\;\;{#2}\;$}\setbox\tw@
 \hbox{$#2$}\ifCD@
 \global\bigaw@\minCDaw@\else\global\bigaw@\minaw@\fi
 \ifdim\wd\z@>\bigaw@\global\bigaw@\wd\z@\fi
 \ifdim\wd\@ne>\bigaw@\global\bigaw@\wd\@ne\fi
 \ifCD@\hskip.5em\fi
 \ifdim\wd\tw@>\z@
 \mathrel{\mathop{\hbox to\bigaw@{\leftarrowfill}}\limits^{#1}_{#2}}\else
 \mathrel{\mathop{\hbox to\bigaw@{\leftarrowfill}}\limits^{#1}}\fi
 \ifCD@\hskip.5em\fi\amp@rs@nd@}
%
%
\newenvironment{CDS}{\@CDS}{\@endCDS}
\def\@CDS{\def\A##1A##2A{\llap{$\vcenter{\hbox
 {$\scriptstyle##1$}}$}\Big\uparrow\rlap{$\vcenter{\hbox{%
$\scriptstyle##2$}}$}&}%
\def\V##1V##2V{\llap{$\vcenter{\hbox
 {$\scriptstyle##1$}}$}\Big\downarrow\rlap{$\vcenter{\hbox{%
$\scriptstyle##2$}}$}&}%
\def\={&\hskip.5em\mathrel
 {\vbox{\hrule width\minCDaw@\vskip3\ex@\hrule width
 \minCDaw@}}\hskip.5em&}
\def\verteq{\Big\Vert&}
\def\novarr{&}
\def\noharr{&&}
\def\SE##1E##2E{\slantedarrow(0,18)(4,-3){##1}{##2}&}
\def\SW##1W##2W{\slantedarrow(24,18)(-4,-3){##1}{##2}&}
\def\NE##1E##2E{\slantedarrow(0,0)(4,3){##1}{##2}&}
\def\NW##1W##2W{\slantedarrow(24,0)(-4,3){##1}{##2}&}
\def\slantedarrow(##1)(##2)##3##4{%
\thinlines\unitlength1pt\lower 6.5pt\hbox{\begin{picture}(24,18)%
\put(##1){\vector(##2){24}}%
\put(0,8){$\scriptstyle##3$}%
\put(20,8){$\scriptstyle##4$}%
\end{picture}}}
\def\vspace##1{\noalign{\vskip##1\relax}}\relax\let\amp@rs@nd@&\iffalse}\fi
 \CD@true\vcenter\bgroup\relax\let\\=\cr\iffalse}\fi\tabskip\z@skip\baselineskip20\ex@
 \lineskip3\ex@\lineskiplimit3\ex@\halign\bgroup
 &\hfill$\m@th##$\hfill\cr}
\def\@endCDS{\cr\egroup\egroup}
%
\newdimen\TriCDarrw@
\newif\ifTriV@
\newenvironment{TriCDV}{\@TriCDV}{\@endTriCD}
\newenvironment{TriCDA}{\@TriCDA}{\@endTriCD}
\def\@TriCDV{\TriV@true\def\TriCDpos@{6}\@TriCD}
\def\@TriCDA{\TriV@false\def\TriCDpos@{10}\@TriCD}
\def\@TriCD#1#2#3#4#5#6{%
\setbox0\hbox{$\ifTriV@#6\else#1\fi$}
\TriCDarrw@=\wd0 \advance\TriCDarrw@ 24pt
\advance\TriCDarrw@ -1em
\def\SE##1E##2E{\slantedarrow(0,18)(2,-3){##1}{##2}&}
\def\SW##1W##2W{\slantedarrow(12,18)(-2,-3){##1}{##2}&}
\def\NE##1E##2E{\slantedarrow(0,0)(2,3){##1}{##2}&}
\def\NW##1W##2W{\slantedarrow(12,0)(-2,3){##1}{##2}&}
\def\slantedarrow(##1)(##2)##3##4{\thinlines\unitlength1pt
\lower 6.5pt\hbox{\begin{picture}(12,18)%
\put(##1){\vector(##2){12}}%
\put(-4,\TriCDpos@){$\scriptstyle##3$}%
\put(12,\TriCDpos@){$\scriptstyle##4$}%
\end{picture}}}
\def\={\mathrel {\vbox{\hrule
   width\TriCDarrw@\vskip3\ex@\hrule width
   \TriCDarrw@}}}
\def\>##1>>{\setbox\z@\hbox{$\scriptstyle
 \;{##1}\;\;$}\global\bigaw@\TriCDarrw@
 \ifdim\wd\z@>\bigaw@\global\bigaw@\wd\z@\fi
 \hskip.5em
 \mathrel{\mathop{\hbox to \TriCDarrw@
{\rightarrowfill}}\limits^{##1}}
 \hskip.5em}
\def\<##1<<{\setbox\z@\hbox{$\scriptstyle
 \;{##1}\;\;$}\global\bigaw@\TriCDarrw@
 \ifdim\wd\z@>\bigaw@\global\bigaw@\wd\z@\fi
 \mathrel{\mathop{\hbox to\bigaw@{\leftarrowfill}}\limits^{##1}}
 }
 \CD@true\vcenter\bgroup\relax\let\\=\cr\iffalse}\fi
 \tabskip\z@skip\baselineskip20\ex@
 \lineskip3\ex@\lineskiplimit3\ex@
 \ifTriV@
 \halign\bgroup
 &\hfill$\m@th##$\hfill\cr
#1&\multispan3\hfill$#2$\hfill&#3\\
&#4&#5\\
&&#6\cr\egroup%
\else
 \halign\bgroup
 &\hfill$\m@th##$\hfill\cr
&&#1\\%
&#2&#3\\
#4&\multispan3\hfill$#5$\hfill&#6\cr\egroup
\fi}
\def\@endTriCD{\egroup}

\title[Algebraic theory of characteristic classes]
{ Algebraic theory of characteristic classes of bundles
with connection} 
\author{ H\'el\`ene Esnault } 
\address{ Universit\"at Essen, FB6 Mathematik, 45 117 Essen, Germany}
\email{ esnault@uni-essen.de}
\thanks{ This work has been partly supported by the DFG Forschergruppe
''Arithmetik und Geometrie''}

\maketitle
\setcounter{section}{-1} 

\section{Introduction}

The Weil algebra homomorphism 
$$
w = \oplus w_n: \oplus_n S^n (\ffg (\C)^*)^{G (\C)}
\to \oplus_n H^{2n}_{DR} (BG (\C))
$$
assigns a de Rham cohomology class of the classifying space
$BG(\C)$ to a $G(\C)$ invariant polynomial on the dual 
$\ffg(\C)^*$ of the Lie algebra associated to 
the $\C$-valued points of an algebraic group $G$ (see
\cite{KN}, Chapter XII). 

In the unpublished note \cite{BK}, Beilinson and Kazhdan give an
algebraic description of $w$ as an iterated Atiyah
extension (see \cite{A}). 

Let $p:\sE \to X$ be a (simplicial) principal $G$ bundle on the
(simplicial) smooth algebraic variety $X$ over a ring $k$ of
characteristic zero. The exact sequence 
$$
0 \to p^* \Omega^{1}_{X} \to \Omega^{1}_{\sE} \to
\Omega^{1}_{\sE/X} \to 0
$$
of regular one forms induces the Atiyah extension 
\begin{equation} \label{At0}
0 \to \Omega^{1}_{X} \to \Omega^{1}_{X, \sE} \to \ffg^{*}_{\sE}
\to 0
\end{equation} 
where 
$$\Omega^{1}_{X, \sE} = (p_* \Omega^{1}_{\sE} )^G,
\ffg^{*}_{\sE} = (p_* \Omega^{1}_{\sE/X})^G = \sE \times_G
\ffg^*.$$ 
Here $G $ acts via the adjoint representation on
$\ffg^*$. 

For example, if $G = GL (r)$, and $E = \sE \times_G k^r$, the
corresponding bundle of $k$-modules $k^r$, then $\ffg^{*}_{\sE} =
{\rm End} \ E$, the endomorphisms of $E$. 

Then \ref{At0} induces a
$n$-extension 
\begin{multline*}
0 \to \Omega^{n}_{X} \to \Lambda^n \Omega^{1}_{X, \sE} \to
\Lambda^{n-1} \Omega^{1}_{X, \sE} \otimes \ffg^{*}_{\sE} \to
\ldots \\
\to \Lambda^{n-i} \Omega^{1}_{X, \sE} \otimes S^i \ffg^{*}_{\sE}
\to \ldots \to S^n \ffg^{*}_{\sE} \to 0
\end{multline*}
which defines a connecting homomorphism
\begin{equation} \label{At}
H^0 (X, S^n \ffg^{*}_{\sE}) \to H^n (X, \Omega^{n}_{X}). 
\end{equation}
Evaluated on $X = BG = (G^{\ell +1} / G)_{\ell}$, the $\Z$-simplicial
scheme classifying $G$-principal bundles, and $\sE = (G^{\ell+1})_{\ell}$, 
the universal $G$-principal bundle, one has 
$$
H^0 (BG, S^n \ffg^{*}_{\sE}) = S^n (\ffg^{*})^G, 
$$
and, if $G$ is reductive, the natural maps
\begin{gather*}
\H^{2n} (BG (L), \Omega^{\geq n}_{BG} ) \to H^n (BG (L),
\Omega^{n}_{BG})\\ 
\H^{2n} (BG (L), \Omega^{\geq n}_{BG} ) \to \H^{2n} (BG (L),
\Omega^{\bullet}_{BG} ) = H^{2n}_{DR}(BG(L))
\end{gather*}
are isomorphisms for any algebraically closed field $L$ of
characteristic 0. 
The Weil homomorphism $w_n$ is just the connecting homomorphism
\ref{At}, where one identifies the right hand side 
with the de Rham cohomology via those two isomorphisms.

Chern-Weil theory assigns to a $C^{\infty}$ manifold $X$ and a
bundle $E$ of rank $r$ with a connection $\nabla$, a morphism 
$$
[\nabla]^* : \oplus_n S^n (\ffg (\C)^*) \to \oplus_n H^0 (X,
\Omega^{2n}_{\infty, {\rm cl}}), 
$$
where $\Omega^{i}_{\infty}$ is the sheaf of $C^{\infty}$ forms
of degree $i$ containing the sheaf $\Omega^{i}_{\infty, {\rm cl}}$ of
closed forms, and $G = GL (r)$, such that 
$$
[\nabla]^* (P) = P (\nabla^2 , \ldots , \nabla^2)
$$
is a closed form, the de Rham class of which is $[E]^*( w(P))$. Here 
$$
[E]^* : \oplus_n H^{2n}_{DR} (BG) \to \oplus_n H^{2n}_{DR} (X)
$$
is the map induced by $E$. 

The theory of secondary classes of Chern-Simons and
Cheeger-Simons is a factorization of $[\nabla]^*$. In a spirit
closed to the algebraic definition of $w$, they both have an
algebraic incarnation. The purpose of this survey is to describe
it. For the analytic side of the theory, we refer to \cite{Edmv}
and \cite{So}.
\section{Chern-Simons theory}  

\subsection{Classical theory} \cite{CS}

Given a $C^{\infty}$ manifold $X$, a bundle $E$ of rank $r$, a
connection $\nabla$, $P \in S^n (\ffg (\C)^*)^{G(\C)}, G = GL
(r)$, Chern and Simons consider the principal $G$ bundle
$$
p : \sE \to X
$$
together with the canonical trivialization $p^* E =
\oplus^{r}_{1} \sO_{\infty}$. In this canonical basis,
$p^*\nabla$ becomes a $r \times r$ matrix $A$ of one forms, with
curvature $F(A)$. They define 
$$
TP(A) = n \int^{1}_{0} P(A, \underbrace{F(tA), \ldots, 
F(tA)}_{(n-1) \ {\rm times}} \ ) dt \in H^0 (\sE, \Omega^{2n-1}_{\infty}), 
$$
a functorial solution to the equation 
$$
d \ ? = P (F(A), \ldots , F(A)). 
$$
In order to define classes of $(E, \nabla)$ (or equivalently of
the local system $E^{\nabla}$) on $X$, and not only on $\sE$,
which of course depends on $E$, they assume $\nabla^2 =0$.
In this case,
$TP(A)$ defines a class $[TP(A)]$ in $H^{2n-1}_{DR} (\sE)$.
Assuming further that $P$ has $\Z$-periods, that is 
$$
P \in {\rm Ker} \ (S^n (\ffg (\C)^*)^{G(\C)} \to H^{2n}_{DR} (BG
(\C), \C/ \Z (n)), 
$$
where $\Z (n) = (2 \pi i)^n \Z$, then the restriction of
$[TP(A)]$ to any $C^{\infty}$ section of $p$ is a well defined
class 
$$
TP (\nabla) \in H^{2n-1} (X, \C /\Z (n)).
$$
We denote it simply by $T c_n (\nabla)$ if $P$ is the polynomial
defining the $n$-th Chern class. 

\subsection{Algebraic Theory} \label{BE}

Let $X$ be a smooth algebraic variety over a field $k$ of
characteristic 0, $(E, \nabla)$ be a bundle of rank $r$ with a
connection, $P \in S^n (\ffg (k)^* )^{G(k)}$, $G = GL(r)$. Given
on a Zariski open set $U \subset X$ a local trivialization of
$E$, $\nabla$ becomes a $r \times r$ matrix of one forms 
$$A \in
H^0 (U, M (r \times r, \Omega^1 )),$$ 
and one can define 
$$
TP(A) \in H^0 (U, \Omega^{2n-1} )
$$
as before with 
$$
d TP(A) = P (\nabla^2, \ldots , \nabla^2 ) \in H^0 (X,
\Omega^{2n}_{{\rm cl}} ).
$$
The point is that if $g \in H^0 (U, GL (r, \sO))$ is a gauge
transformation, then 
$$
TP(A) - TP (g^{-1} dg + g^{-1}  Ag)
$$
is locally exact for $n \geq 2$, and thus viewed as a class in
$H^0 (X, \Omega^{2n-1} / d \Omega^{2n-2})$, or in the subgroup
$H^0 (X, \sH^{2n-1}_{DR}) \subset H^0 (X, \Omega^{2n-1} / d
\Omega^{2n-2})$ if $\nabla^2 =0$, with
$\sH^{2n-1}_{DR} = \Omega^{2n-1}_{{\rm cl}} / d \Omega^{2n-2}$, it
depends only on $\nabla$. We denote it by 
$$
TP (\nabla) \in H^0 (X, \Omega^{2n-1} / d \Omega^{2n-2} )
$$
for $n\geq 2$. (In \cite{BE}, it is denoted by $w_n(E, \nabla,
P)$). If $P$ is the polynomial
defining the $n$-th Chern class, we simply denote it by
$T c_n (\nabla)$. . 

Due to the Bloch-Ogus theory \cite{BO}, the restriction maps to
the generic point are injective: 
$$
H^0 (X, \Omega^{2n-1} / d \Omega^{2n-2}) \subset
\Omega^{2n-1}_{k (X)} / d \Omega^{2n-2}_{k (X)}
$$
$$
H^0 (X, \sH^{2n-1}_{DR} ) \subset H^{2n-1}_{DR} (k (X)), 
$$
for $n \geq 2$. 

Thus $TP( \nabla)$ is recognized at the generic point of $X$.
One asks what kind of algebraic class of $E$ it controls. 

To this aim, one first shows that $TP(\nabla)$ in fact lies in a
subgroup $$E^{0,2n-1}_{n} \subset H^0 (X, \Omega^{2n-1} / d
\Omega^{2n-2})$$ arising from the coniveau spectral sequence.
When $k$ is algebraically closed, the group $E^{0,2n-1}_{n}$
maps to the group 
$$
CH^{n}_{1, DR} (X) = \left( \frac{\{ \oplus \Z Z, Z \ \
\mbox{prime cycle of codimension} \ \ n \}}{\sim_{DR}} \right)
\otimes_{\Z} k 
$$
where $Z \sim_{DR} 0$ if there is a divisor $W$ containing the
support of $Z$ such that the class of $Z$ 
in $H^{2n}_{DR,W} (X)$, the de Rham cohomology of $X$ with
supports along $W$, vanishes.
Clearly,
$CH^{n}_{1,DR} (X)$ is a quotient of the Chow group with 
$k$-coefficients $CH^n (X) \otimes_{\Z} k$, and for $n =2$, it
coincides with the Griffiths group $\otimes k$. 

We denote by $c_P (E)$ the class of $E$ associated to $P$ in
$CH^n (X) \otimes_{\Z} k$. 

\begin{prop} (see \cite{BE}, Proposition 5.4.1) Let $k$ be an
algebraically closed field of characteristic 0. 
The image of $TP(\nabla)$ in $CH^{n}_{1,
DR} (X)$ equals the image of $c_P (E)$. 
\end{prop}

[Strictly speaking, the proof in loc. cit. assumes $k = \C$,
$\nabla^2 =0$, and deals with 
$$
CH^{n}_{1, {\rm Betti}} (X) = \frac{\{ \oplus \Z Z, Z \ \
\mbox{prime cycle of codimension} \ \ n \}}{\sim_{\rm Betti}}
$$
where $Z \sim_{\rm Betti} 0$ if there is a divisor $W$
containing the support of $Z$ such that the class of $Z$
in $H^{2n}_{{\rm Betti}, W} (X, \Z
(n))$, the Betti cohomology of $X$ with supports
along $W$, vanishes. Then $CH^{n}_{1, {\rm Betti}} (X)$ is a
quotient of $CH^n (X)$. 

It is straightforward to generalize. First, the sequence
(5.4.3) of loc. cit. has an obvious de Rham version. Then one has
that 
\begin{multline*}
Tc_n (\nabla) \in {\rm Im} \ \H^{n} (X, \sK_n \> d \ {\rm log} >>
\Omega^n \to \ldots \to \Omega^{2n-1}) \\
\> d >> H^0 (X, \Omega^{2n-1}
/ d \Omega^{2n-2})
\end{multline*}
(\cite{Eadc}, section 2.2). The map $d$ factors through
$$
\H^n (X, 0 \to 0 \to \Omega^{n+1} / d \Omega^n \to \Omega^{n+2}
\to \ldots \to \Omega^{2n-1})
$$
and 
\begin{multline*}
E^{0,2n-1}_n = {\rm Im} \ \H^n (X, 0 \to 0 \to \Omega^{n+1} / d
\Omega^{n} \to \Omega^{n+2} 
\to \ldots \to \Omega^{2n-1} )\\
\to
H^0 (X, \Omega^{2n-1} / d \Omega^{2n-2}).
\end{multline*}
Finally, to prove
that the image of $TP(\nabla)$ in $CH^{n}_{1, DR} (X)$ is
the correct one, one replaces loc. cit. (5.4.5) by the exact
sequence 
\begin{multline*}
0 \to (0 \to \Omega^n / d \Omega^{n-1} \to \Omega^{n+1} \to
\ldots \to \Omega^{2n-2}_{{\rm cl}} ) 
\\ 
\to (\sK_n \to \Omega^n /d \Omega^{n-1} \to \Omega^{n+1} \to
\ldots \to \Omega^{2n-1} ) 
\\
\to \sK_n \oplus \Omega^{2n-1} / d \Omega^{2n-2} [-n] \to 0.
\end{multline*}
Here and in section 2, $\sK_n$ is the Zariski sheaf 
$$
{\rm Im} \ (\sK^{M}_{n} \to i_{k (X)^*} K^{M}_{n} (k (X)))
$$
where $i_{k(X)}: {\rm Spec} \ k (X) \to X$ is the inclusion of
the generic point and $K^{M}_{n}$ is the Milnor $K$-theory.]

The main theorem is now 

\begin{thm} (see \cite{BE}, Theorem 5.6.2) \label{rig}
Let $X$ be projective smooth over $\C$,
$$
P \in {\rm Ker} \ (S^n (\ffg (\C)^*)^{G(\C)} \to H^{2n} (BG(\C),
\C/\R (n))), 
$$
$n \geq 2$, and let $(E, \nabla)$ be a flat bundle on $X$. Then
$TP(\nabla) =0$ if and only if the image of $c_P (E)$ in
$CH^{n}_{1, {\rm Betti}} (X) \otimes_{\Z} \R$ vanishes. 
\end{thm} 

\subsection{Question}
When $\nabla$ is not flat, there are many examples of 
non-vanishing 
$TP (\nabla)$ classes. However, when $\nabla$ is
flat, we don't know any $(n \geq 2)$. 

This raises the question of whether $TP(\nabla)$ always vanishes
under the assumption of Theorem \ref{rig}. It is related to a
question arising from Nori's work (\cite{N}, Introduction),
often called Nori's conjecture: if 
$$Z\in {\rm Ker} \ (CH^2 (X) \to
H^{4}_{\sD} (X, \Q (2))),$$ 
where $X$ is projective smooth over
$\C$, and $H^{a}_{\sD} (b)$ is the Deligne cohomology, is $Z
\otimes \Q$ algebraically equivalent to 0?

Since for $n = 2$,
$\sim_{\rm Betti}$  is the algebraic
equivalence, a generalization of Nori's question to any
codimension $n$ is: if 
$$Z \in {\rm Ker} \ (CH^n (X) \to
H^{2n}_{\sD} (X, \Q (n))),$$ 
where $X$ is projective smooth over
$\C$, does the image of $Z$ in $CH^{n}_{1, {\rm Betti}} (X)
\otimes \Q$ vanish? 

Reznikov's theorem \cite{R} asserts that if $X$ is projective smooth
over $\C$ and $E$ is flat on $X$, then 
$$
c_n (E) \in {\rm Ker} \ (CH^n (X) \to H^{2n}_{\sD} (X, \Q (n)))
$$
for $n \geq 2$. In view of this and of \ref{rig}, the question
for $Z= c_n (E)$, $E$ flat, is equivalent to the question of the
vanishing of $T c_n (\nabla)$ for any flat structure $\nabla$ on
$E$. 

\section{Cheeger-Simons theory}

\subsection{Classical theory: differential characters}
\cite{CSI} 

Chern-Weil theory provides an invariant $P(\nabla^2, \ldots,
\nabla^2) \in H^0 (X, \Omega^{2n}_{\infty, {\rm cl}} )$ when $\nabla^2
\neq 0$ and Chern-Simons theory an invariant $TP(\nabla) \in
H^{2n-1} (X, \C / \Z (n))$ when $\nabla^2 =0$ and $P$ has 
$\Z$-periods. Chern-Simons theory defines invariants for $(E,
\nabla)$ combining the two. To this aim, Chern and Simons
define the ancestor of Deligne cohomology, the group of
differential characters 
$$
\hat{H}^{2n} (X) = \H^{2n} (X, \Z (n) \to \sO_{\infty} \to
\ldots \to \Omega^{2n-1}_{\infty}), 
$$
which is an extension of 
$$
{\rm Ker} \ (H^0 (X, \Omega^{2n}_{\infty, {\rm cl}} ) \to H^{2n} (X,
\C / \Z (n)))
$$
with $H^{2n-1} (X, \C / \Z (n))$. (The notation differs from
theirs, as well as the hypercohomology presentation). 

This group is functorial, and there is a natural ring structure
on
$
\oplus_n \hat{H}^{2n} (X).
$
For $P$ with $\Z$-periods, they define $\hat{c}_P \in
\hat{H}^{2n} (X)$, by saying that in bounded rank and dimension,
there is a classifying space for bundles with connection. Since
a $C^{\infty}$ bundle always carries a $C^{\infty}$ connection,
this space also classifies bundles and thus has no odd
dimensional cohomology. This implies that on this space, the
differential characters inject into closed differential forms.
Chern-Weil theory provides then the classes of the universal
connection. 

\subsection{Algebraic theory: algebraic differential characters}
\cite{Eadc}

Let $X$ be a smooth algebraic variety defined over a field $k$
of characteristic 0. One defines the group of algebraic
differential characters by 
$$
AD^n (X) : = \H^n (X, \sK_n \> d \ {\rm log} >> \Omega^n \to
\ldots \to \Omega^{2n-1} ), 
$$
an extension of 
$$
{\rm Ker} \ (H^0 (X, \Omega^{2n}_{{\rm cl}} ) \to \frac{\H^{2n} (X,
\Omega^{\geq n})}{{\rm Im} \ CH^n (X)} )
$$
by 
$$
\H^n (X, \Omega^{\infty} \sK_n ): = \H^n (X, \sK_n \> d \ {\rm
log} >> \Omega^n \to \ldots \to \Omega^{\dim X} ). 
$$
This group is functorial, and $AD (X) = \oplus_n AD^n (X)$ has a
natural ring structure. Over $k = \C$, it maps to $\hat{H}^{2n}
(X)$, compatibly with the extension. In general, it also maps to
$CH^n (X)$ and $H^0 (X, \Omega^{2n-1} / d \Omega^{2n-2} )$ for
$n \geq 2$. For $n =1$, $AD^1 (X)$ is the group of isomorphism
classes of $(E, \nabla)$, where $E$ is a rank 1 bundle and
$\nabla$ is a connection, and $AD^1 (X) \to H^0 (X,
\Omega^{2}_{{\rm cl}})$ is the curvature map. 

\begin{thm} \label{adc} \cite{Eadc}
There are functorial classes $c_n (E, \nabla) \in AD^n (X)$,
such that $c_1 (E, \nabla)$ is the class of $({\rm det} \ E,
{\rm det} \ \nabla)$ in $AD^1 (X)$, such that $c_n (E, \nabla)$
lifts $c_n (E) \in CH^n (X)$, $Tc_n (\nabla ) \in H^0 (X,
\Omega^{2n-1} / d \Omega^{2n-2} )$ for $n \geq 2$, and
$\hat{c}_n (E, \nabla) \in \hat{H}^{2n} (X)$ if $k= \C$. Those
classes verify the Whitney product formula. 
\end{thm} 

Classically, there are two ways of constructing classes of
vector bundles in a cohomology theory: via the splitting
principle, knowing $c_1 (\sO (1))$, where $\sO (1)$ is the
tautological bundle on the projective bundle $\pi: \P = \P (E)
\to X$ to $E$, and knowing the freeness of the cohomology of
$\P$ as a module over the cohomology of $X$. Or via the
(simplicial) classifying space $BG$, and defining classes of the
universal (simplicial) bundle $EG$. For connections, $AD(\P)$ is
not a free module over $AD(X)$, $\sO (1)$ does not have a
connection, nor does $EG$ have an algebraic connection. 

Nonetheless, one can define the classes $c_n (E, \nabla)$ via a
modified splitting principle and via a universal $BG$
construction as well. 

For the splitting principle, one has to consider connections
with values in a differential graded algebra. This leads to a
definition of general groups of algebraic differential
characters fulfilling the splitting principle, and thereby the
unicity of the classes of Theorem \ref{rig}. As an illustration,
let us explain the objects for $n=2$. $\nabla$ defines a
splitting 
$$
\tau : \Omega^{1}_{\P} \to \pi^* \Omega^{1}_{X}
$$
of the exact sequence of forms such that $\tau \circ \nabla$ is
compatible with $\pi^* E \to \sO (1)$. Thus it defines a class
$$
\xi = (\sO (1), \tau \circ \nabla) \in \H^1 (\P, \sK_1 \> \tau
\circ d \ {\rm log} >> \pi^* \Omega^{1}_{X}). 
$$
One observes that with respect to the splitting $\tau$, one has 
$$
d (\Omega^{2}_{\P /X} ) \subset \Omega^{3}_{\P /X} \oplus
\Omega^{2}_{\P /X} \otimes \pi^* \Omega^{1}_{X} \oplus
\Omega^{1}_{\P /X} \otimes \pi^* \Omega^{2}_{X} . 
$$
In particular, the complex
$$
\sK_2 \> d \ {\rm log} >> \frac{\Omega^{2}_{\P}}{\Omega^{2}_{\P
/X}} \> d >> \pi^* \Omega^{3}_{X}
$$
is well defined. With respect to a natural product $\cup$, its
cohomology equals $AD^2 (X) \oplus AD^1 (X) \cup \xi$ (see
\cite{Eadc}, section 2). 

For the universal construction, one
first constructs a cohomology associated to a smooth
(simplicial) algebraic variety $X$ and a (simplicial) vector
bundle $E$, which depends on $E$, such that its value on $(BG,
EG)$ is ``tautological'', and such that a connection $\nabla$ in
$E$ maps this cohomology to $AD^n (X)$. 

In the analytic context, Beilinson und Kazhdan (\cite{BK})
developed this procedure to recover the Cheeger-Simons classes
in $\hat{H}^{2n} (X)$. They construct an analytic cohomology
associated to $X$ and $E$, which maps to $\hat{H}^{2n} (X)$ via
a connection in $\nabla$ on $E$. The main ingredient, which they
construct and which we use in our algebraization of their
construction, is the following filtered differential graded
algebra. 

Let $p: \sE \to X$ be the principal $G$ bundle to $E$, and
$\Omega^{1}_{X, \sE}$ be as in \ref{At}. Then 
\begin{gather*}
\Omega^{n}_{X,\sE} = \oplus_{a+b=n} \Omega^{a,b}_{X, \sE} 
\\
\Omega^{a,b}_{X, \sE} = \Lambda^{a-b} \Omega^{1}_{X, \sE}
\otimes S^b \ffg^{*}_{\sE} 
\\
F^n \Omega^{\bullet}_{X, \sE} = \oplus_{a \geq n}
\Omega^{a,b}_{X, \sE} [a+b]. 
\end{gather*}
There is a natural differential $\Omega^{n}_{X, \sE} \to
\Omega^{n+1}_{X, \sE}$ extending the K\"ahler differential, and
the natural injection 
$$
(\Omega^{\bullet}_{X}, \Omega^{\geq n}_{X}) \to
(\Omega^{\bullet}_{X, \sE} , F^n \Omega^{\bullet}_{X, \sE})
$$
is a filtered quasi-isomorphism. A connection $\nabla$ is
equivalent to an inverse quasi-isomorphism 
$$
[\nabla ] : \Omega^{\bullet}_{X, \sE} \to \Omega^{\bullet}_{X}
$$
and the integrability condition is equivalent to 
$$
[\nabla ] (F^n \Omega^{\bullet}_{X, \sE}) \subset \Omega^{\geq
n}_{X} .
$$
The Weil homomorphism can be understood as a map 
$$
S^n (\ffg^*)^G \> w_n >> F^n \Omega^{\bullet}_{X, \sE} [2n]
$$ 
and the tautological cohomology needed is 
$$
\H^n (X, {\rm cone} \ ( \sK_n \otimes S^n (\ffg^*)^G [-n] \> d
\ {\rm log} \ \oplus - w_n >> F^n \Omega^{\bullet}_{X, \sE} [n])
[-1])
$$
(\cite{Eadc}, section 3). 

\subsection{Questions}

\begin{enumerate}
\item[1.] Since there are flat bundles $(E, \nabla)$ such that $c_n (E)
\neq 0$ in $CH^n (X) \otimes_{\Z} \Q$, for example $E =
\oplus^r L_i $, $L_i \in {\rm Pic}^0 X, X$ abelian variety
(\cite{B}), there are classes $c_n (E, \nabla)$ which are not
vanishing for flat bundles. The question is what is the coniveau
of those classes, that is the largest codimension $a$ such that
a codimension $a$ subvariety $Z \subset X$ exists with $c_n (E,
\nabla) | X-Z =0$. This question is related to the vanishing of
$Tc_n (\nabla)$ in section 1. 
\item[2.] The group $AD^n (X)$ mixes $\sK$-cohomology with
cohomology of differential forms. It would be more powerful to
mix $\sK$-cohomology with something related to Betti or \'etale
cohomology. This probably would solve the question on the
vanishing of $Tc_n (\nabla)$. 
\end{enumerate}

\section{Riemann-Roch theorems}

Let $f: X \to S$ be a projective smooth morphism of relative
dimension $d$ over an algebraically closed field $k$ of
characteristic 0, with $S$ and $X$ smooth. Mumford \cite{M}
observed that if $d =1$, then 
\begin{equation} \label{van}
c_n (R^1 f_* (\Omega^{\bullet}_{X/S} )) = 0 \ \ \mbox{in} \ \
CH^n (S) \otimes_{\Z} \Q, 
\end{equation}
applying the Grothendieck-Riemann-Roch theorem to the single
sheaves of the relative de Rham complex
$\Omega^{\bullet}_{X/S}$. The same argument shows 
\begin{multline} \label{c}
c_n (\sum_i (-1)^i R^i f_* (\Omega^{\bullet}_{X/S} \otimes E,
\nabla_{X/S})) = \\
(-1)^d f_* (c_d (\Omega^{1}_{X/S} ) \cdot c_n
(E)) \ \mbox{in} \ CH^n (S) \otimes_{\Z} \Q ,
\end{multline}
if $E$ is a bundle on $X$, endowed with a flat relative
connection $\nabla_{X/S}$. If $\nabla_{X/S}$ comes from a flat global
connection $\nabla$, the Gau\ss-Manin
bundles $$R^i f_* (\Omega^{\bullet}_{X/S} \otimes E)$$ carry the
Gau\ss-Manin connection $GM (\nabla)$. 

If $k = \C$, the work of Bismut-Lott \cite{BL} and Bismut
\cite{Bi} shows that \ref{c} is true in $\hat{H}^{2n} (X)
\otimes \Q$: 
\begin{multline} \label{cs} 
\hat{c}_n (\sum_i (-1)^i [R^i f_* (\Omega^{\bullet}_{X/S} \otimes
E, \nabla_{X/S}), GM (\nabla)]) = \\
(-1)^d f_* (c_d
(\Omega^{1}_{X/S} ) \cdot \hat{c}_n (E, \nabla) ) \ \mbox{in} \
\hat{H}^{2n} (X) \otimes \Q . 
\end{multline}

For $n=1$, there is an analogy with the situation where $S$ is a
finite field $\F_q$, $(E, \nabla)$ is a local system $V$. Then
the work of Deligne \cite{D}, \cite{DI}, and subsequent work by
Laumon \cite{L}, S. Saito \cite{SS}, T. Saito \cite{TS}, show
that \ref{cs} for $n=1$ remains true: 
\begin{equation} \label{det}
{\rm det} \ \sum_i (-1)^i H^{i}_{{\rm \acute{e}t}} (X, V) = (-1)^d {\rm
det} \ V | c_d (\Omega^{1}_{X})
\end{equation}
as Frobenius-modules over $\F_q$. 

Both \ref{cs} and \ref{det} involve classes of the local system.
The classes $Tc_n (\nabla)$ and $c_n (E, \nabla)$ reflect also
the choice of the algebraic structure $E$. One shows 

\begin{thm} \cite{BEII} \label{rr1}
\begin{enumerate}
\item[a)] 
\begin{multline*}
c_1 (\sum_i (-1)^i [R^i f_* (\Omega^{\bullet}_{X/S} \otimes E ,
\nabla_{X/S}), GM (\nabla)]) = \\
(-1)^d f_* (c_d(\Omega^{1}_{X/S} )
\cdot c_1 (E, \nabla)) \ \mbox{ in } \ AD^1 (S) \otimes \Q 
\end{multline*}
\item[b)] 
\begin{multline*}
Tc_n (\sum_i (-1)^i [R^i f_* (\Omega^{\bullet}_{X/S}\otimes E ,
\nabla_{X/S}), GM (\nabla) ]) = \\
(-1)^d f_* (c_d (\Omega^{1}_{X/S} ) \cdot Tc_n (E, \nabla))
\ \mbox{ in } \ H^0 (S, \sH^{2n-1}_{DR}) \ \mbox{ if } \ n \geq 2. 
\end{multline*}
\end{enumerate}
\end{thm}

More generally, Theorem \ref{rr1} remains true if one replaces $\nabla$
by a flat connection with logarithmic poles along a relative
normal crossing divisor $Y = \bigcup_i Y_i$, $c_d
(\Omega^{1}_{X/S})$ by the relative top Chern class 
$$
c_d (\Omega^{1}_{X/S} ({\rm log} \ Y), {\rm res}_{Y_i} ) \in
\H^d (X, \sK_d \to \oplus_i \sK_d |Y_i \to \oplus_{i <j} \sK_d |
Y_i \cap Y_j \to \ldots )
$$
as defined by T. Saito in \cite{TS}, 
using the existence of the residue maps
${\rm res}_{Y_i}: \Omega^1_{X/S}(\log Y) \to \sO_{Y_i}$,
and the classes $c_n (E,
\nabla), Tc_n (\nabla)$ by the corresponding classes involving
logartihmic poles, which we haven't discussed here at all.

The introduction of logarithmic poles is necessary in order to
have sufficiently many bundles with connections with the help of
which one
can reduce the problem to curves as in \cite{L}, \cite{SS},
\cite{TS}. 

Another generalization of Theorem \ref{rr1} is 
\begin{thm} \cite{BEII} \label{rr2}
Under the assumptions 
\begin{enumerate}
\item[(i)] $S = {\rm Spec} \ K$, $K$ function
field over $k$ and 
\item[(ii)] $\nabla^2 \in \H^0 (X, f^*
\Omega^{2}_{S} \otimes {\rm End} \ E)$, 
\end{enumerate}
one has the same
conclusion as in Theorem \ref{rr1}. 
\end{thm} 

Note that the assumption (ii) of \ref{rr2} allows to define
Gau\ss-Manin connections on $S$. Under the weaker assumption
$\nabla^{2}_{X/S} =0$, no natural connection is defined on $R^i
f_* (\Omega^{\bullet}_{X/S} \otimes E, \nabla_{X/S} )$. However
one has 
\begin{thm} \cite{BEII} \label{rr3}
Under the assumptions \ref{rr2} (i) and $d=1$ 
(thus (ii) is automatically fulfilled), 
there is a naturally defined connection on 
$$ {\rm det} (\sum_{i=0}^{i=2}
(-1)^i R^i f_* (\Omega^{\bullet}_{X/S} \otimes E,
\nabla_{X/S})).$$
With respect to this connection, the
conclusion of Theorem \ref{rr1}, a) holds true. 
\end{thm}

In $AD^n (S) \otimes_{\Z} \Q$, one obtains the Riemann-Roch
formula in the very trivial case $X = Y \times S$,
$f=$ projection, by deforming the relative de Rham cohomology to
some relative Higgs cohomology, and keeping the Gau\ss-Manin
connection \cite{Err}. 

\subsection{Remarks and Questions} 
\begin{enumerate}
\item[1)] Mumford's observation \ref{van} generalizes to 
$$
c_n (R^1 f_* \Omega^{\bullet}_{X/S} ) =0 \ \ \mbox{in} \ \ CH^n
(S) \otimes_{\Z} \Q, 
$$
as computed by van der Geer \cite{vdG}, applying the
Grothendieck-Riemann-Roch theorem to powers of a principal
polarization of the corresponding family of abelian varieties.
Then \ref{c} implies 
$$
c_n (R^2 f_* \Omega^{\bullet}_{X/S}) =0 \ \ \mbox{in} \ \ CH^n
(S) \otimes_{\Z} \Q
$$
for a family of smooth projective surfaces as well. Does one have 
$$
c_n (R^i f_* \Omega^{\bullet}_{X/S} ) =0 \ \ \mbox{in} \ \ CH^n
(S) \otimes_{\Z} \Q
$$
for a smooth projective family over a smooth base $S$?

The contribution of the singularities of $f$ at infinity is
difficult to understand. Mumford \cite{M} shows actually that
the Gau\ss-Manin bundle with logarithmic singularities
of a semi-stable curve
has torsion Chern classes in the Chow group. In general, even the analytic
statement in Deligne cohomology is not understood (see
\cite{Edmv}, 3.6 Questions). 
\item[2)] Deligne \cite{DI} gave a proof of \ref{det} for $d=1$,
rank $V=1$ using the geometry of the abelian variety ${\rm
Pic}^0 X$. In view of the shape of the Riemann-Roch formula,
containing the expression $c_d (\Omega^{1}_{X/S})$, it would be
natural to try to understand it using the geometry of the moduli
of Higgs bundles, Hitchin's map and Higgs cohomology.
\item[3)] Gau\ss-Manin bundles, when $\nabla^2 =0$, are direct
images of special holonomic $\sD$-modules with regular
singularities under special projective morphisms. A general
Riemann-Roch formula would require a Chern-Simons and
Cheeger-Simons theory for $\sD$-modules, the algebraic cycle
part of which should be given by \cite{LI}. 
\end{enumerate}

\subsection{Acknowledgements:}

B. Ang\'eniol taught me the meaning of Atiyah classes
(\cite{EI}, Introduction), before he left mathematics. I'd like
to thank him again for the discussions we have had at that time.
A large part of the ideas exposed loosely in this survey go back
to joint work with S. Bloch. It is my pleasure to thank him for
his generosity and the joy of the mathematical exchanges. My
first understanding of de Rham complexes goes back to joint work
with E. Viehweg. It influenced later my work on classes. I thank
him for his constant support. \bibliographystyle{plain}
\renewcommand\refname{References}
 
\end{document}